\documentclass[12pt]{amsart}
\usepackage{amssymb,amsmath}
\makeatletter
\def\@cite#1#2{{\m@th\upshape\bfseries%
[{#1\if@tempswa{\m@th\upshape\mdseries, #2}\fi}]}}
\makeatother
\theoremstyle{plain}
\newtheorem{thm}{Theorem}[section]
\newtheorem{lem}[thm]{Lemma}
\newtheorem{cor}[thm]{Corollary}

\theoremstyle{definition}
\newtheorem{rem}[thm]{Remark}

\newtheorem{note}[thm]{Note}
\newtheorem{eg}[thm]{Example}
\newtheorem{ques}[thm]{Question}
\newtheorem{ans}[thm]{Answer}

\newcommand{\Prf}{\noindent\textbf{Proof.\ }}
\newcommand{\bx}{\strut\hfill$\blacksquare$\medbreak}

\newcommand{\ca}{\mathrm{C}^*}

\newcommand{\ol}{\overline}
\newenvironment{sbmatrix}{\left[\begin{smallmatrix}}{\end{smallmatrix}\right]}
\newenvironment{spmatrix}{\left(\begin{smallmatrix}}{\end{smallmatrix}\right)}

\newcommand{\wot}{\textsc{wot}}

\newcommand{\bbB}{{\mathbb{B}}}
\newcommand{\bbC}{{\mathbb{C}}}

\newcommand{\bbR}{{\mathbb{R}}}
\newcommand{\bbT}{{\mathbb{T}}}

 \newcommand{\B}{{\mathcal{B}}}

\renewcommand{\H}{{\mathcal{H}}}

 \newcommand{\K}{{\mathcal{K}}}

 \newcommand{\M}{{\mathcal{M}}}
 \newcommand{\N}{{\mathcal{N}}}
\renewcommand{\O}{{\mathcal{O}}}

\renewcommand{\phi}{\varphi}
\newcommand{\upchi}{{\raise.35ex\hbox{$\chi$}}}
\newcommand{\fA}{{\mathfrak{A}}}

\newcommand{\fS}{{\mathfrak{S}}}

\newcommand{\qand}{\quad\text{and}\quad}

\newcommand{\qfor}{\quad\text{for}\quad}

\newcommand{\qiff}{\quad\text{if and only if}\quad}
\newcommand{\ran}{\operatorname{Ran}}

\newcommand{\spn}{\operatorname{span}}

\newcommand{\sumin}{\sum_{i=1}^n}
\newcommand{\ai}{A_i}

\newcommand{\aistar}{A_i^*}

\newcommand{\cp}{\Phi}
\newcommand{\bofk}{\B(\K)}
\newcommand{\bofh}{\B(\H)}
\newcommand{\rowa}{(A_1, \ldots, A_n)}
\newcommand{\rows}{(S_1, \ldots, S_n)}

\newcommand{\cpmap}{completely positive map}
\newcommand{\kj}{\K_j}
\newcommand{\xjk}{X_{jk}}
\newcommand{\wjk}{W_{jk}}
\newcommand{\ehij}{A_{i,j}}
\newcommand{\ehik}{A_{i,k}}

\begin{document}

\title[Completely Positive Maps]%
{Quantum channels, wavelets, dilations, and representations of $\O_n$}
%
\author[D.W.Kribs]{David~W.~Kribs${}^{1}$}
\thanks{{\it key words and phrases. } Hilbert space, operator, completely
positive map, quantum channel, orthogonal wavelet, isometric dilation,
representation, Cuntz algebra}
\thanks{${}^{1}$partially supported by a Canadian NSERC Post-doctoral
Fellowship.}
\address{Department of Mathematics and Statistics, University of
Guelph, Guelph, ON, CANADA N1G 2W1} \email{dkribs@uoguelph.ca}

\date{}
\begin{abstract}
We show that the representations of the Cuntz $\ca$-algebras $\O_n$ which
arise in wavelet analysis and dilation theory can be classified through a
simple analysis of completely positive maps on finite-dimensional
space. Based on this
analysis, we find an application in quantum information theory; namely, a
structure theorem for the
fixed point set of a unital quantum channel. We also include some open
problems motivated by this work.
\end{abstract}
\maketitle

There has been considerable recent interest in the analysis of completely
positive maps on finite-dimensional space. There are a number of reasons
for this including  connections with
wavelet analysis \cite{BJwave,BEJ,Jorgen_min}, dilation theory
\cite{DKS,JK},  representation theory of the Cuntz $\ca$-algebras $\O_n$
\cite{BJwave,BEJ,DKP,DKS}, and quantum information theory
\cite{Werner,King,Ruskai,Shor}. The results
obtained in
the current paper have implications for each of these areas.
In presenting this work, another goal we have is to push further the
connections between the various areas mentioned above.

In the first section we establish a result for completely positive
maps. While we focus on the finite-dimensional setting, this is not
necessary in the proof. A structure theorem for the fixed point set of a
unital quantum channel is contained in the second section.
In particular, we prove that the fixed point set is a $\ca$-algebra which
is equal to the commutant of the algebra generated by any choice of row
contraction which determines the channel.
We discuss the 2-dimensional channels \cite{King,Ruskai}, and use the
theorem
to classify them by their fixed point sets.
The representation theory for $\O_n$ is considered in the third
section. We focus on a subclass of representations arising in
dilation theory and wavelet analysis
\cite{BJwave,BEJ,DKP,DKS,Jorgen_min}. Each of these
representations determines a completely positive map on finite-dimensional
space. We ask
if these representations can be classified just by examining the map. An
affirmative answer is provided by the result on completely positive maps
from the first section. Finally, in the fourth section we pose some open
questions which are motivated by this work.

\section{Completely Positive Maps} \label{S:main}

In this section we present a theorem for completely positive maps on
finite-dimensional space.
Let $\K$ be a finite-dimensional Hilbert space and let $\bofk$ be the
bounded operators on $\K$. It is well-known (see
\cite{Choi,Kraus,Pa_cb_book} for instance)
that every completely positive map $\Phi : \bofk \rightarrow \bofk$
is determined by a
row matrix $A = \rowa$ of operators in $\bofk$ in the following sense:
\begin{eqnarray}
\Phi (X) &=& \sumin \ai X \ai^*  \qfor X\in\bofk.
\end{eqnarray}
The map is {\it unital} if, in addition, $\Phi(I) = \sum_{i=1}^n A_i A_i^*
= AA^* =I$. Of course,
there will be many choices of tuples $A$ which determine
a given completely
positive map in this way.  We prove the following result.

\begin{thm}\label{main}
Let $\K$ be a finite-dimensional Hilbert space. Let $\Phi : \bofk
\rightarrow \bofk$ be a
completely positive unital map, and suppose $A = \rowa$ is a choice  of
operators which determine $\Phi$ as in $(1)$. If $p$ is an orthogonal
projection in $\bofk$, then we have the following
equivalences for the range space $\ran(p)$ of the projection:
\begin{itemize}
\item[$(i)$] $\Phi (p) \geq p$ if and only if $\ran (p)$ is $\ai^*$-invariant for $1 \leq i \leq n$.
\item[$(ii)$] $\Phi (p) \leq p$ if and only if $\ran (p)$ is $\ai$-invariant for $1 \leq i \leq n$.
\item[$(iii)$] $\Phi (p) = p$ if and only if $\ran (p)$ is $\ai$-reducing for $1 \leq i \leq n$.
\end{itemize}
Furthermore, for a given projection $p$, one of these equivalences holds
for a particular choice of $A$ in $(1)$ if
and only if it holds for all choices of $A$ in $(1)$.

In each of these cases the limit
\[
\Phi^\infty (p) := \lim_{k\rightarrow\infty} \Phi^k (p)
\]
exists. This operator belongs to the fixed point set $\bofk^\Phi$
of $\Phi$, and it is computed from the fixed point set as follows.
\begin{itemize}
\item[$(i)'$] If $\Phi (p) \geq p$, then $\Phi^\infty (p) = \inf \{
X\in\bofk^\Phi : X\geq p\}.$
\item[$(ii)'$] If $\Phi (p) \leq p$, then $\Phi^\infty (p) = \sup \{
X\in\bofk^\Phi : X\leq p\}.$
\item[$(iii)'$] If $\Phi (p) = p $, then the infimum and supremum from
$(i)'$ and $(ii)'$
are both equal to $p$.
\end{itemize}

\end{thm}

While the existence of the limit $\Phi^\infty (p)$ and the final three
facts seem to be of independent interest, the applications of
this result contained in the rest of the paper are based on $(i)-(iii)$.
\begin{note}
After preparing this article the author discovered that conditions
$(ii)$ and $(iii)$ were recently established in \cite{BJKW}, although the
proofs here are different.
\end{note}
Before proving the theorem, we establish a lemma which generalizes a
result from \cite{DKS}.

\begin{lem} \label{espaces}
Let $X$ be a positive operator in $\bofk$ which satisfies the inequality
$0 \leq X \leq \cp (X).$
Then the eigenspace $\ker (X -  ||X||I )$ is $\aistar$-invariant for $1
\leq i \leq n$.
\end{lem}

\Prf
Without loss of generality assume that $||X|| = 1$. Let $\M = \ker (X - I)$. Then for any vector
$\xi \in \M$,
\begin{eqnarray*}
  \| \xi \|^2  = ( X \xi, \xi) &\leq& \sum_{i=1}^n ( A_i X A_i^* \xi, \xi ) \\
    & = & \sum_{i=1}^n ( X A_i^* \xi, A_i^* \xi )
     \leq \sum_{i=1}^n ( A_i^* \xi, A_i^* \xi ) = \| \xi \|^2.
\end{eqnarray*}
In particular, all the inner product inequalities are actually equalities. Since $X$ is a positive
contraction, the
only way this can happen is if each of the vectors $\aistar \xi$ belongs to $\M$. Hence $\M$
is $\aistar$-invariant.
\bx

{\noindent}{\bf Proof of Theorem~\ref{main}.}
First note that the equivalences $(i)$ and $(ii)$ are duals of each other. Indeed, since $\Phi$ is
unital,
\[
\Phi (p) \geq p \qiff \Phi (I - p) \leq I - p,
\]
and the subspace $\ran (p)$ is $\ai^*$-invariant precisely when $\ran (I-p) = \ran(p)^\perp$ is
$\ai$-invariant. Thus we shall prove $(i)$.

If $\Phi (p) \geq p$, then an application of the previous lemma
yields the $\ai^*$-invariance of $\ran (p)$. To see the converse
of $(i)$, with respect to the orthogonal decomposition $\K = p\K
\oplus p^\perp\K$, suppose $\ai$ can be written in matrix form as
\[
\ai = \left[
\begin{matrix}
 B_i & 0 \\ C_i & D_i
\end{matrix} \right]
\qfor 1 \leq i \leq n.
\]
Then the relation $\Phi (I_\K) = \sumin \ai \aistar = I_\K$ yields the identities:
\[
\sumin B_i B_i^* = I_{p\K}, \!\!\! \quad\text{}\!\! \quad \sumin B_i C_i^* =0
\qand \sumin (C_iC_i^* + D_i D_i^*) = I_{p^\perp\K}.
\]
Thus upon writing $p = \begin{sbmatrix} I & 0 \\ 0 & 0 \end{sbmatrix}$ with respect to this
spatial decomposition we get
\[
\cp(p) = \sumin A_i p A_i^* = \left[
\begin{matrix}
 I & 0 \\ 0 & \sumin C_i C_i^*
\end{matrix} \right]
\geq  \left[
\begin{matrix}
 I & 0 \\ 0 & 0
\end{matrix} \right]    = p.
\]

For $(iii)$, notice that the previous computation shows that $\Phi(p) = p$ when $\ran(p) = p\K$
is  $\ai$-reducing. Indeed,
the subspace $p\K$ is reducing for the operators  $A= ( A_1, \ldots, A_n )$  exactly
when $C_i = 0$ for $1\leq i \leq n$. On the other hand, if $\Phi (p) =p$,
then the inequalities in $(i)$ and
$(ii)$ are satisfied. Hence $\ran(p)$ is $\ai$-reducing.

Concerning the limits, suppose $p$ satisfies equivalence $(i)$. Since $\Phi$ is positive and
unital, we have
\[
0 \leq p \leq \Phi(p) \leq \Phi^2 (p) \leq \ldots \leq \Phi^k (p) \leq \ldots \leq I.
\]
This is a monotone increasing sequence of positive operators which is bounded above, hence we
obtain the existence of the limit $\Phi^\infty (p) =
\lim_{k\rightarrow\infty} \Phi^k (p)$. It is clear from the form $(1)$
that $\Phi$ is continuous, hence $\Phi^\infty (p)$ is
fixed under the action of $\Phi$. On the other hand, when $X \geq p$
is fixed by the map, we have  $X \geq \Phi^k(p)$ for $k\geq 1$, showing
that $\Phi^\infty (p)$ is bounded above by
every fixed point which majorizes $p$.  It follows that  $\Phi^\infty (p)$
actually {\it is} the infimum.

The proof of $(ii)'$ is analogous since, in that case, the operators
$\Phi^k(p)$ form a decreasing
sequence of positive operators which are majorized by $p$.  Finally, when $\Phi (p) =p$, we get
$\Phi^\infty (p) =p$ so that
$(i)'$ and $(ii)'$ show that the infimum and supremum are both equal to
$p$.
\bx


\begin{rem}
While our focus is on the finite-dimensional case, we note that
the proof of Theorem 1.1 works for {\it any} completely positive
unital map determined as in $(1)$ by a row contraction $A$. The
only change in the conclusion is that the limit $\Phi^\infty (p)$
converges in the strong operator topology.
\end{rem}

\section{Quantum Channels}\label{S:quantum}

Mathematically, a  {\it quantum channel} is a
completely positive trace preserving map on finite-dimensional space.
In the language of quantum information theory, a channel describes the
transfer of quantum information, or {\it qubits}, from `Alice' to `Bob'
(see
\cite{Werner,King,Ruskai,Shor} for
some recent related analysis, as well as the text \cite{CN} for
general information). The operators $A = \rowa$ which determine a quantum
channel $\Phi$ as in $(1)$ are called the {\it Kraus operators}
\cite{Kraus} of the channel.

In general the fixed point set $\bofk^\Phi= \{ X\in\bofk | \Phi(X) =
X\}$ of a completely
positive map is just  a self-adjoint subspace. In particular, generally
the fixed
point set is {\it not} closed under multiplication. We obtain
the following
structure theorem for the fixed point set of a unital quantum channel.

\begin{thm}\label{quantum}
Let $\Phi : \bofk \rightarrow \bofk$ be a unital quantum channel.
If $A = \rowa$ determines $\Phi$ as in $(1)$, then the algebra $\fA$
generated by $A_1, \ldots, A_n$ is a $\ast$-algebra which depends only on
$\Phi$. Further, the fixed point set $\bofk^\Phi$ of the map $\Phi$
coincides with the commutant of $\fA$,
$\fA^\prime = \{ X\in\bofk \,\,|\,\, A_iX=XA_i \qfor 1 \leq i \leq n\}$,
and hence
is itself a
$\ast$-algebra containing the identity operator on $\K$.
\end{thm}

We begin by pointing out a special case of the theorem, established
previously in \cite{BJKW,DKS}, which holds more generally. For
completeness we provide a proof.

\begin{lem}\label{irred}
Let $\K$ be a finite-dimensional space and let $\Phi: \bofk \rightarrow
\bofk$ be a unital completely positive map which is determined as in
$(1)$ by $A=\rowa$. If $\fA = \bofk$, then the fixed point set for $\Phi$
consists of scalars,
$\bofk^\Phi = \bbC I$.
\end{lem}

\Prf
Suppose $X=X^*$ is non-scalar and satisfies $\Phi(X)=X$. Then Lemma 1.3
can
be adapted to show that the eigenspaces corresponding to the two extremal
eigenvalues for $X$ are perpendicular $\aistar$-invariant subspaces, which
are both non-trivial since $X$ is non-scalar. Thus the algebra $\fA^*$ has
proper invariant subspaces, a contradiction.
\bx

The key observation for unital quantum channels is that
invariant subspaces for the determining $n$-tuples are actually reducing.

\begin{lem}\label{reducing}
Let $\Phi : \bofk \rightarrow \bofk$ be a unital quantum channel. Then
every projection which satisfies
$\Phi (p) \geq p$,  also satisfies $\Phi (p) =p$. Thus if $A=\rowa$
determines $\Phi$, then every subspace which is invariant for the family
$\{ A_i^* : 1 \leq i \leq n \}$ is also reducing for the family.
Furthermore, every subspace which is invariant for $\{A_i : 1 \leq i \leq
n\}$ is also reducing.
\end{lem}

\Prf
Since $\Phi$ is positive and unital, we have $0 \leq \Phi(p) \leq I$. Thus if we are given a
projection $p$ which satisfies $p \leq \Phi (p)$, then trace preservation
ensures  we must have
equality, $p = \Phi(p)$. A similar analysis follows for the dual notion
$\Phi(p) \leq
p\leq I$. The rest of the lemma follows from Theorem 1.1.
\bx

{\noindent}{\bf Proof of Theorem~\ref{quantum}.} We first show
$\fA$ is a $\ast$-algebra. Let $\{p_j\}$ be a maximal family of
pairwise orthogonal projections which are each minimal reducing
for the family $\{A_i : 1 \leq i \leq n\}$. Let $\kj=p_j\K$ for
each $j$. Then $A_i p_j = p_j A_i$ for all $i,j$, and by
maximality $I_\K =\sum_j p_j$. Hence the algebra $\fA = \sum_j p_j
\fA p_j$ is block diagonal with respect to this family, and the
blocks $\fA p_j = p_j \fA = p_j \fA p_j$ are algebras themselves.
But the trivial subspaces are the only subspaces which are
invariant for $\fA p_j$. Indeed, if a subspace $p\K$ is invariant
for each $A_i$ and $p\leq p_j$ is supported on some $p_j$, then
$p\K$ is reducing for the $A_i$ by Lemma~\ref{reducing}, and hence
by the minimality of $p_j$ we have either $p=0$ or $p=p_j$. It
follows from Burnside's classical theorem that the restricted
finite-dimensional algebra $\fA p_j$ is equal to  $\B(\kj)$, in
particular the algebras $\fA p_j$ are self-adjoint. Therefore
$\fA=\fA^*$ is a $\ast$-algebra. (In particular, it is a
finite-dimensional $\ca$-algebra.)

Since $\Phi$ is unital, that is $\Phi(I) =I$, it is clear that every $X$
which
commutes with $A_1, \ldots ,A_n$ will be fixed by $\Phi$. Hence the fixed
point set
$\bofk^\Phi$ contains  the commutant $\fA^\prime$.

To see the converse, let $\ehij = A_i p_j$ for all $j$ and
$i=1,\ldots, n$, and for each $j$ put $B_j = (A_{1,j},\ldots
,A_{n,j})$. Given $X\in\bofk$, let $X=(\xjk)$ where $\xjk= p_j X
p_k$. If $X$ satisfies $\Phi(X)= X$, then a computation shows that
$\Phi_j(X_{jj}) = X_{jj}$ where $\Phi_j :
\B(\kj)\rightarrow\B(\kj)$ is the unital completely positive map
$\Phi_j(Y) = \sumin \ehij Y \ehij^*$. However, as observed above,
the minimality of $p_j$ as an $\ai$-invariant subspace gives us
$\fA p_j = \B(\kj)$. Hence by Lemma~\ref{irred}, we have $X_{jj} =
x_{jj}p_j$ for some scalar $x_{jj}$.

For $j\neq k$, we claim that either $\xjk=0= X_{kj}$ for all
$X=X^*$ in $\bofk^\Phi$,
or there is a unitary $\wjk : \K_k\rightarrow \K_j$ with
$\ehij = \wjk \ehik \wjk^*$ for $i= 1,\ldots ,n$. Thus suppose there is an
$X=X^*$ with
$\Phi(X) =X$ and $\xjk\neq 0$, normalized so that $||\xjk|| = 1$. Let
$\M = \{ \xi\in\K_k : ||\xjk\xi|| = ||\xi|| \}$, and let $\N = \xjk\M$ be
the corresponding subspace of $\kj$. Then for $\xi\in\M$ we have
\[
\xjk\xi = (p_j \Phi(X) p_k )\xi = \Phi(\xjk)\xi = (B_j \xjk^{(n)} B_k^*)
\xi.
\]
This implies that each $\ehik^*$ leaves $\M$ invariant. Indeed,
since $B_j$ and $B_k^*$ are contractions, and $\xjk$ achieves its
norm on $\xi$, it follows that $B_k^*\xi$ belongs to the subspace
$\M^{(n)}$ on which $\xjk^{(n)}$ achieves its norm. Thus by Lemma
2.3, $\M$ is an $\ai$-reducing subspace contained in $\K_k$. Since
$\M$ is non-zero, the minimality of $p_k$ gives $\M=\K_k$. By
considering $X_{kj} = \xjk^*$, we see that $\N = \kj$ also. In
particular, $\xjk$ and $\xjk^*$ are partial isometries, and it
follows that $\wjk = \xjk|_{\K_k}: \K_k \rightarrow\kj$ is a
unitary operator.

The identity above shows that $\wjk = B_j \wjk^{(n)} B_k^*$. Hence for
$\xi\in\K_k$
\[
||\xi|| = ||\wjk \xi|| = ||B_j\wjk^{(n)}B_k^* \xi|| \leq ||\wjk^{(n)}
B_k^* \xi|| \leq ||\xi||.
\]
Thus $B_j$ acts as an isometry from the range $\ran \wjk^{(n)}B_k^*$ onto
the range $\ran\wjk = \kj$. As $B_j$ is a row contraction, it must be zero
on the orthogonal complement of $\ran \wjk^{(n)}B_k^*$. Hence, $B_j^* $ is
an isometry from $\kj$ onto $\ran \wjk^{(n)} B_k^*$. Consequently,
$B_j^* \wjk = \wjk^{(n)}B_k^*$, in other words, $\ehij^* = \wjk \ehik^*
\wjk^*$ for $i= 1, \ldots ,n$, as claimed.

Now suppose $\Phi(Y)=Y = Y^* = (Y_{jk})$. If $j\neq k$ is a pair which
fits into the above analysis, then
\begin{eqnarray*}
\wjk \Phi_k (\wjk^*Y_{jk}) &=& \sumin \wjk \ehik \wjk^* Y_{jk} \ehik^* \\
&=& \sumin \ehij Y_{jk} \ehik^* \\
&=& p_j \Phi(Y) p_k = Y_{jk}.
\end{eqnarray*}
Thus by Lemma~\ref{irred}, $\wjk^*Y_{jk}$ is scalar. Hence $Y_{jk}
= y_{jk} \wjk$ for some scalar $y_{jk}$, and also $Y_{kj} =
Y_{jk}^* = \ol{y}_{jk} \wjk^*$. The other off-diagonal entries of
$Y$ are either zero, or have a similar form.

This analysis gives us a handle on the matrix entries in the decomposition
$X=(\xjk)$ for $X=X^* $ in $\bofk^\Phi$. The corresponding form for each
$A_i$ in this decomposition is given by $A_i = \sum_j A_ip_j = \sum_j
\ehij$, and hence a computation shows that $XA_i =A_iX$ for $i=1,\ldots
,n$ and
$X=X^* $ in $\bofk^\Phi$. Since the self-adjoint subspace $\bofk^\Phi $ is
spanned by its self-adjoint part, it follows that $\fA^\prime$ contains
the entire fixed point set $\bofk^\Phi$. This completes the proof.
\bx

\begin{note}
We note that our approach in proving Theorem 2.1 was motivated by work in
\cite{DKS} on dilation theory, which is discussed in the next
section. This connection leads to an open question posed in Section
4 which, if answered,  could mesh
the theory of quantum channels with that of certain operator algebras on
infinite-dimensional space.
\end{note}

We finish this section by showing how this theorem can be used to classify
quantum channels by their fixed point sets.

\begin{eg}
The quantum channels on $\M_2 = \B(\bbC^2)$, the so called {\it qubit
channels}, were recently characterized in
\cite{Ruskai}. The identity matrix together with the {\it Pauli matrices}
$\{ I, \sigma_x, \sigma_y, \sigma_z \}$ form a basis for $\M_2$ where
\[
\sigma_x = \begin{spmatrix} 0 & 1 \\ 1 & 0 \end{spmatrix} \,\,\,\,\,\,\,\,
\sigma_y = \begin{spmatrix} 0 & -i \\ i & 0 \end{spmatrix}
\,\,\,\,\,\,\,\,
\sigma_z = \begin{spmatrix} 1 & 0 \\ 0 & -1 \end{spmatrix}.
\]
Every unital qubit channel is
equivalent, through unitary conjugations at both the input and output
stages, to a map $\Phi$ which is diagonal and real with
respect to this basis. In other words; $\Phi(I)=I$, $\Phi(\sigma_x) =
\lambda_1 \sigma_x$, $\Phi(\sigma_y) = \lambda_2 \sigma_y$,
$\Phi(\sigma_z) = \lambda_3 \sigma_z$, where $\lambda_1, \lambda_2,
\lambda_3 \in \bbR$ (for more on the geometry of channels see the text
\cite{CN}).

Thus describing the set of channels on $\M_2$ amounts to deriving
conditions
on the $\lambda_k$ which guarantee such a diagonal map is completely
positive and trace preserving. The norm of a completely positive map
$\Phi$ is given by $||\Phi|| = || \Phi(I) ||$. Hence for these
diagonal maps, a simple necessary condition for complete positivity is
that each $|\lambda_k |\leq 1$. Using Choi's lemma \cite{Choi}, necessary
and sufficient conditions
on the $\lambda_k$ were computed in \cite{Ruskai} which lead to a
description of the entire
set of qubit channels, though the conditions are rather technical.

From Theorem 2.1 and the basic theory of finite-dimensional
$\ca$-algebras,  we know there are just three possibilities for the fixed
point algebra $\M_2^\Phi$ of a unital qubit channel $\Phi$. Somewhat
surprisingly, we don't need the extra
conditions on $\{\lambda_1, \lambda_2, \lambda_3\}$ from \cite{Ruskai} to
classify the fixed point sets.
\end{eg}

\begin{cor}
Let $\Phi : \M_2 \rightarrow \M_2$ be a unital quantum channel.
Then the fixed point set for $\Phi$ is a $\ast$-subalgebra of
$\M_2$ containing the identity operator and satisfying one of the
following conditions.
\begin{itemize}
\item[$(i)$] $\M_2^\Phi = \M_2$ if and only if each $\lambda_k = 1$
and $\Phi$ is the identity map.
\item[$(ii)$] $\M_2^\Phi = \bbC I$ if and only if $\lambda_k \neq 1$ for
$k= 1,2,3$.
\item[$(iii)$] $\M_2^\Phi = \spn \{ |v_1\!\!>\!<\!\!v_1|,
|v_2\!\!>\!<\!\!v_2| \}$, where
$\{ v_1, v_2 \}$ is an orthonormal basis for $\bbC^2$ and
$|v_i\!\!>\!<\!\!v_i|$
is the rank
one projection of $\bbC^2$ onto $\spn \{ v_i \}$, and this case holds if
and only if exactly one of $\{\lambda_1,\lambda_2,\lambda_3\}$ is equal to
$1$.
\end{itemize}
\end{cor}

\Prf
By the previous discussion we may assume $\Phi$ is diagonal with respect
to the Pauli basis and that $|\lambda_k| \leq 1$ for $k=1,2,3$. An
elementary
computation shows that the Pauli basis has the property that projections
$p=p^2=p^* \in \M_2$ are either trivial ($p=0$ or $p=I$) or represented as
\[
p = \frac{1}{2} I + a \sigma_x + b \sigma_y + c \sigma_z,
\]
where $a,b,c\in\bbR$ and $a^2 + b^2 + c^2 = 1 / 4$. Thus $\Phi(p) =p$
exactly when $a=a\lambda_1$, $b=b\lambda_2$ and $c=c\lambda_3$.

In particular, the case $\lambda_1 = \lambda_2 = \lambda_3=1$
corresponds to $a,b,c$ being free variables and $\Phi$ being the
identity map on $\M_2$. By Theorem 2.1 the operator algebraic
characterizations of $\M_2^\Phi$ in $(ii)$ and $(iii)$ are the
only remaining possibilities for the fixed point set. Further, the
finite-dimensional $\ca$-algebra  $\M_2^\Phi$ is spanned by its
projections. Thus $\M_2^\Phi = \bbC I$ exactly when there are no
non-trivial projections fixed by $\Phi$. Equivalently, there are
no solutions $a,b,c$ to the above identities. Clearly this holds
if and only if each $\lambda_k \neq 1$.

Hence the remaining cases must satisfy
$(iii)$, and we claim this is when exactly one of the $\lambda_k$ is
equal to 1.
Indeed, either one or two of the $\lambda_k$ must be equal to 1 for
the case $(iii)$ maps, since otherwise we would be in one of the first two
cases by the previous paragraph. If, say, $\lambda_1=\lambda_2=1$, then
$\Phi$ fixes $\sigma_x$ and $\sigma_y$. Whence, $(-i)\sigma_x \sigma_y =
\sigma_z$ is also fixed by $\Phi$ since the fixed point set is an algebra.
But this would imply that $\lambda_3 = 1$, and we are really in case
$(i)$. Similarly, it is easy to see that when any two of
$\{\lambda_1,\lambda_2,\lambda_3\}$ are equal to $1$, the third must be
as
well. Thus the operator algebra characterization of $\M_2^\Phi$ in $(iii)$
occurs precisely when there is exactly one $\lambda_k$ equal to 1.
\bx

\section{Applications To Representation Theory For $\O_n$}
\label{S:background}

In this section we show that  representations of the Cuntz
$\ca$-algebra $\O_n$ arising from
dilations and wavelets can be classified through an analysis of completely
positive maps.
Given a positive integer $n \geq 2$, $\O_n$ is the universal
$\ca$-algebra generated by the relations
\[
 s_i^*  s_j = \delta_{ij} I \qfor 1 \leq i,j \leq n \qand \sum_{i=1}^n s_i
s_i^* =I.
\]
An $n$-tuple $S = \rows$ of operators in $\bofh$ which satisfies these relations consists of
isometries with pairwise orthogonal ranges, for which the range projections of the isometries
span the entire (necessarily infinite-dimensional) Hilbert space $\H$. Up
to isomorphism,
$\O_n$ is the $\ca$-algebra generated by {\it any} such $n$-tuple since it is simple. A theorem
of Glimm's \cite{Glimm} suggests it is not possible to find a meaningful classification of all
representations of $\O_n$ (it is an `NGCR' algebra). However, there are
good reasons for
attempting to classify particular subclasses of these representations, including connections with
the study of endomorphisms of $\bofh$, finitely correlated states, dilation theory, wavelet
analysis, and the theory of non-selfadjoint operator algebras (See
\cite{BJwave,BEJ,DKP,DKS,Jorgen_min} for examples from  different
perspectives).

The representations $\pi$ of $\O_n$ on a space $\H$ which are of interest
here have the property that there exists a finite-dimensional subspace
$\K$ of $\H$ which is co-invariant and cyclic for the isometries $S =
\rows$, where $\pi(s_i) = S_i$. In other words,
\begin{itemize}
\item[$(i)$] $S_i^* \K \subseteq \K \qfor 1\leq i \leq n$
\item[and]
\item[$(ii)$] $\H = \bigvee_{i_1,\ldots ,i_k} S_{i_1}\cdots
S_{i_k} \K$,
where the closed span is over all indices $1\leq i_1, \ldots, i_k \leq n$
and $k \geq 1$.
\end{itemize}

We let  $A = \rowa$ be the row contraction of matrices consisting of the
compressions to $\K$ of the isometries, so that
$A_i = P_\K S_i|_\K = (S_i^*|_\K)^*$. Notice that $\sum_{i=1}^n A_i A_i^* =
I_\K$ and thus  $\Phi (X) = \sum_{i=1}^n A_i X A_i^*$ defines a completely
positive unital map on $\bofk$.
These representations were classified in \cite{DKS}. They form
the subclass of representations of $\O_n$ which arise through the
{\it minimal isometric dilations} of row contractions of matrices
\cite{Bun,Fra1,Pop_diln}. They also include the representations of $\O_n$
which come from wavelet analysis.

Every {\it orthogonal wavelet} of scale $n$ is determined by a
scaling function $\phi$ in $ L^2 (\bbR)$ that determines functions
which generate a `wavelet basis' for  $L^2 (\bbR)$. On the other
hand, the Fourier expansion of $\phi$ also determines so called
wavelet filter functions $m_1, \ldots, m_n$ in  $L^\infty (\bbT)$.
Let $\rho$ be a primitive $n$th root of unity. The orthogonality
of the wavelet is embedded in the statement that the complex
matrices $\frac{1}{\sqrt{n}} \big( m_i (\rho^k z) \big)_{i,k=1}^n$
are unitary for  ${\rm a.a.}\,\,z\in \bbT$. Given such a wavelet,
a representation of $\O_n$ on $L^2 (\bbT)$ is obtained by defining
isometries $S_i f (z) = m_i (z)\, f(z^n)$ for $1\leq i \leq n$.

Extensive analysis has been conducted on these and other related wavelet
representations (see \cite{BJwave,BEJ,Jorgen_min,JK} for instance).
When the scaling function $\phi$ is compactly supported, the associated
representation possesses a finite-dimensional subspace $\K$
which satisfies $(i)$ and $(ii)$ above for the isometries $S_i$. In fact,
these representations are quite specialized in that $\K$
can be chosen to be spanned by Fourier basis vectors.
Thus, the orthogonal wavelet representations form a {\it subclass}
of the $\O_n$ representations arising through dilation theory.
Hence, the analysis in \cite{DKS} can be applied to these
representations, and it has been recently by Jorgensen
\cite{Jorgen_min}.

The decomposition theory of \cite{DKS} can
quickly become computationally cumbersome.
This basic problem provided the initial motivation for this paper:

\begin{ques}\label{classify}
Given a row contraction $A = \rowa$ of matrices, or equivalently a
completely
positive map $\Phi$ on finite-dimensional space, is it possible to
classify the associated representation of $\O_n$ just in terms of  $\Phi$,
without reference to the $n$ operators $\ai$?
\end{ques}

It was shown in \cite{DKS} that decomposing these representations, which
act  on infinite-dimensional space, amounts to an exercise in
finite-dimensional matrix algebra. This
matrix algebra essentially consists of identifying a maximal family of
pairwise orthogonal minimal $\ai^*$-invariant subspaces.
Each of these minimal `anchor' subspaces generates an irreducible subspace
for the $S_i$, and with this information the representations can
be classified. Thus we have the following:

\begin{ans}
Theorem~\ref{main} answers Question~\ref{classify} in the {\it
affirmative.} In particular, such anchor subspaces are identified
by finding any maximal family of mutually orthogonal minimal projections
satisfying $p \leq \Phi(p)$, and these subspaces {\it only} depend on
$\Phi$. Thus all the classification results for these representations in
previous papers \cite{BJwave,BEJ,DKS,Jorgen_min} can be restated without
reference to the minimal
$\ai^*$-invariant subspaces, only to the minimal projections satisfying
this inequality. For instance, this gives a new
characterization of irreducibility: one of these
representations is irreducible precisely when there is a {\it unique} minimal
projection $p$ satisfying $p \leq \Phi(p)$.
\end{ans}

We show how this new perspective can ease the computational burden
by considering an example from each of the dilation and wavelet
settings.

\begin{eg}
The following  example is due to Arveson, and  appeared in the
seminal paper \cite{ArvII} as an example of a completely positive
map for which the fixed point set  $\bofk^\Phi$ is {\it not} an
algebra. Nonetheless, it provides a satisfying application of the
method introduced here. For $k \geq 2$, let $\cp : \M_k
\rightarrow \M_k$ be the completely positive unital map defined by
(assuming an orthonormal basis for $\bbC^k$ has been fixed)
\[
\cp ( [ x_{ij} ]) =
\begin{sbmatrix}
 x_{11} & & &   0 \\  & \ddots  & &  \\
& & x_{k-1k-1} & \\
0 & &  & \frac{1}{k-1} \sum_{i=1}^{k-1} x_{ii}
\end{sbmatrix} .
\]

Without the new perspective discussed above, classifying a representation
of $\O_n$ generated by $\Phi$ through dilation theory would first require
finding a row contraction $A = \rowa$
which determines $\Phi$ as in $(1)$. Next, the minimal $\aistar$-invariant
subspaces would have to be computed. With the new perspective this
becomes a triviality: The rank one projections
$E_{1,1}, \ldots, E_{k-1,k-1}$ clearly satisfy $E_{ii} \leq \cp
(E_{ii})$, and are obviously
minimal with respect to this property. This is all the information we need
here to describe a representation of $\O_n$ generated by $\Phi$,
and we emphasize that {\it no} reference is required to a row contraction
$A$ which determines $\Phi$.

Indeed, it follows that such a representation breaks up into the
direct sum of $k-1$ irreducible subrepresentations. In particular, the
representation will be irreducible if and only if $k=2$. The ranges
of $E_{1,1}, \ldots, E_{k-1,k-1}$ provide one-dimensional `anchor'
subspaces which generate the irreducible subspaces associated with the
irreducible subrepresentations. Further, some thought shows
that the free semigroup algebra $\fS$ (the $\wot$-closed algebra generated
by the isometries from the dilation \cite{DKP,DKS,DP1}) associated
with each subrepresentation is unitarily equivalent to the tractable
`one-dimensional atomic' free semigroup algebra arising in the literature
\cite{DP1}. Thus the free
semigroup algebra of the full representation is unitarily equivalent
to the direct sum of $k-1$ copies of this algebra.
\end{eg}

\begin{eg}
Matrix representations were worked out in \cite{BJwave} for the {\cpmap}s $\Phi$
determined by wavelet representations of $\O_3$ with $\K = \spn \{ z^0,
z^{-1}, z^{-2} \} \subseteq L^2
(\bbT)$. The authors show how
combining an eigenvalue analysis of the matrix, together with a
computation of the fixed point set for the map, can be used to discern
information on the representation.
The method presented in this section stream-lines this analysis in that
it can be used to obtain this information in one fell swoop.

Consider  the ordered basis of standard matrix units for $\M_3$
\[
\bbB = \{ E_{0,0}, E_{-1,-1}, E_{-2,-2}, E_{0,-1}, E_{0,-2}, E_{-1,-2}, E_{-2,-1}, E_{-1,0},
E_{-2,0} \},
\]
corresponding to the ordered basis
$\{  z^0, z^{-1}, z^{-2} \}$ of $\K$. An example of a matrix representation
$[ \Phi ]_{\bbB}$ in this basis for a {\cpmap} $\Phi$ determined by a wavelet representation of
$\O_3$ is given by
\[
[ \Phi ]_{\bbB} =
\left[
\begin{matrix}
1 & 0& 0& 0& 0& 0& 0& 0&  0 \\
0 & \frac{1}{\sqrt{2}} & 1- \frac{1}{\sqrt{2}}& 0& 0& 0& 0& 0&  0 \\
0 & \frac{1}{\sqrt{2}} & 1-\frac{1}{\sqrt{2}}& 0& 0& 0& 0& 0&  0 \\
0 & 0 & 0& \frac{1}{\sqrt{2}} & 0& -\frac{1}{\sqrt{2}}& 0& 0& 0 \\
0 & 0 & 0& \frac{1}{\sqrt{2}} & 0& -\frac{1}{\sqrt{2}}& 0& 0& 0 \\
0 &   \frac{1}{\sqrt{2}} & -\frac{1}{\sqrt{2}} & 0& 0 & 0 & 0& 0& 0 \\
0 &   \frac{1}{\sqrt{2}} & -\frac{1}{\sqrt{2}} & 0& 0 & 0 & 0& 0& 0 \\
0 & 0 &0 & 0& 0 & 0 & 0& 0& 0 \\
0 & 0 &0 & 0& 0 & 0 & 0& 0& 0
\end{matrix}
\right]
\]
(for the reader of \cite{BJwave}, this example is generated by
taking
$g =2$, $N=3$, $\lambda_0 =1$, $\lambda_1 = \frac{1}{\sqrt{2}} = \lambda_2$).

The rank one projection $E_{0,0}$ satisfies $\Phi (E_{0,0}) = E_{0,0}$, and  is clearly minimal
with this property.
Solving for projections
$p = (p_{ij})$, with ranges orthogonal to that of $E_{0,0}$ (so $p_{ij} =  0$ if $i=0$ or $j=0$),
and such that $p \leq \Phi(p)$,  yields the inequalities
\[
p_{11} \leq \frac{p_{11}}{\sqrt{2}} + p_{22}(1 - \frac{1}{\sqrt{2}})  \qand
p_{22} \leq \frac{p_{11}}{\sqrt{2}} + p_{22}(1 - \frac{1}{\sqrt{2}}).
\]
Hence $p_{22} \geq p_{11}$ and $p_{11} \geq p_{22}$, so that equality is
achieved. Further elementary analysis shows that $p_{12}=p_{21}=0$. Thus
$\cp ( E_{-1,-1} + E_{-2,-2}) = E_{-1,-1} + E_{-2,-2}$ and this rank
two projection is minimal satisfying $p\leq \Phi(p)$ since
any smaller projection satisfying this inequality
would have to also satisfy the above inequalities involving $p_{11}$ and $p_{22}$.

It now follows that the associated wavelet representation has two
irreducible subrepresentations, and the corresponding irreducible
subspaces are generated by the anchor subspaces
\[
(E_{-1,-1} + E_{-2,-2})\K \qand E_{0,0}\K .
\]
Also, $\O_3$ cyclic vectors for the corresponding irreducible summands
can be obtained simply by taking bases for the generating anchor
subspaces. In particular, $\{ z^0 \}$ and $\{ z^{-1}, z^{-2} \}$ will
suffice for the two subspaces.
\end{eg}



\section{Open Questions}\label{S:open}

We finish by taking the opportunity to pose some open problems motivated
by the work in this paper.

\begin{ques}
As Remark 1.4 points out, the equivalences $(i)-(iii)$ in Theorem 1.1 are
valid when the $A_i$ act on infinite-dimensional space. Thus we ask
whether there are subclasses of representations of $\O_n$ arising from
dilations of {\it infinite} rank $n$-tuples $A$ which can be classified
using just the associated completely positive map $\Phi$? At present there
is not even a non-trivial subclass of such representations which has been
classified in any way.

A natural class to consider could be wavelet representations of $\O_n$ for
which the associated scaling function is {\it not} compactly
supported. It can be shown that these representations  have a
co-invariant cyclic subspace
$\K$ which is  infinite-dimensional, hence they do indeed arise
from the dilations of infinite rank $n$-tuples. These representations also
have the advantage of
having an explicit formula for the generating isometries, which
is determined by the wavelet filter functions and the scaling function.

We mention that this class could provide examples that shed light on a
deep problem in free semigroup algebra theory, posed in \cite{DKP}, which
is related to the invariant subspace problem.
\end{ques}

\begin{ques}
In many respects, quantum information theory is still in its infancy. For
instance, there are certainly connections with operator theory and
operator algebras, but these seem to be underdeveloped at present.

A
natural question to ask here is whether the representations of $\O_n$
determined by quantum channels through dilation theory have a meaningful
interpretation in quantum information theory?
Could they provide a theory for describing the `external noise'
associated with  quantum transmissions?
\end{ques}

\begin{ques}
We also wonder how the distinguished fixed points $\Phi^\infty (p)$ from
Theorem 1.1 fit into the analysis of completely positive maps? A
connection with eigenvalue analysis is suggested by Lemma 1.3, but we are
unable to say anything substantial at this point.
\end{ques}


{\it Acknowledgements.} The author is grateful to Palle Jorgensen
and Stephen Power for helpful discussions. Thanks also to members
of the Department of Mathematics at the University of Iowa, and
Department of Mathematics and Statistics at Lancaster University
for kind hospitality during the preparation of this article.

\vspace{0.1in}

Note Added in Proof. We mention that Theorem~2.1 has motivated
forthcoming work \cite{HKL} on quantum error correction.
Furthermore, a new simpler proof of Theorem~2.1 is included in
\cite{HKL}.

\end{document}